\documentclass[12pt]{article}

\newcommand{\es}{\mathop{\rm ess \; sup}\limits}
\newcommand{\ei}{\mathop{\rm ess \; inf}\limits}

\def\dual{\,^{^{\complement}}\!}

\def\Rn{{\mathbb{R}^n}}
\def\a{\alpha}

\def\o{\omega}

\usepackage{color}
\usepackage{srcltx}
\usepackage{amsmath}
\usepackage{latexsym}
\usepackage{amsfonts}
\usepackage{amssymb}
\usepackage{times}

\usepackage{latexsym}
\newcounter{theorem}

\renewcommand{\thetheorem}{\arabic{section}.\arabic{theorem}}

\newenvironment{thm}[1]{\par\addvspace{0.5cm}
    \begin{sloppypar}\refstepcounter{theorem}%
    {\bf #1 \thetheorem.}\it{}}{\end{sloppypar}}

\newcommand{\eh}{\hfill}\newlength{\sperr}

\newenvironment{theorem}{\begin{thm}{Theorem}} {\end{thm}}
\newenvironment{lemma}{\begin{thm}{Lemma}} {\end{thm}}
\newenvironment{corollary}{\begin{thm}{Corollary}} {\end{thm}}

\newenvironment{defi}[1]{\par\addvspace{0.5cm}
\begin{sloppypar}\refstepcounter{theorem}%
{\bf #1 \thetheorem.}\rm{}}{\end{sloppypar}}

\newenvironment{definition}{\begin{defi}{Definition}}{\end{defi}}

\newenvironment{remark}{\begin{defi}{Remark}}{\end{defi}}

\newenvironment{proof}{{\settowidth{\sperr}{\rm Proof}
\par\addvspace{0.3cm}\parbox[t]{1.3\sperr}{\rm P\eh r\eh o\eh o\eh f\eh. }%
}}{\nopagebreak\mbox{}\hfill $\Box$\par\addvspace{0.25cm}}

\newcommand{\al}{\alpha}
\newcommand{\bt}{\beta}
\newcommand{\dl}{\delta}
\newcommand{\om}{\omega}
\newcommand{\Om}{\Omega}
\newcommand{\lb}{\lambda}

\newcommand{\ve}{\varepsilon}
\newcommand{\gm}{\gamma}

\newcommand{\vi}{\varphi}

\newcommand{\M}{\mathcal{M}}

\newcommand{\intl}{\int\limits}

\newcommand{\vs}{\vspace}

\newcommand{\diam}{\mbox{\,\rm diam\,}}

\renewcommand{\thetheorem}{\arabic{section}.\arabic{theorem}}

\begin{document}

\begin{center}
\LARGE Maximal, potential and singular operators\\ in the
local "complementary" variable exponent Morrey type spaces
\end{center}

\centerline{\large {\bf Vagif S. Guliyev
}, \,\,  {\bf Javanshir J. Hasanov}
}

\centerline{\it Institute of Mathematics and Mechanics, Baku,
Azerbaijan.}

\centerline{\it E-mail: vagif@guliyev.com}

\centerline{\it E-mail: hasanovjavanshir@yahoo.com.tr}

\

\centerline{\large \bf Stefan G. Samko}

\centerline{\it University of Algarve, Portugal.}

\centerline{\it E-mail: ssamko@ualg.pt}

\

\begin{abstract}
We consider local "complementary" generalized Morrey spaces
${\dual \cal M}_{\{x_0\}}^{p(\cdot),\om}(\Om)$ in which the $p$-means of   function are controlled over $\Om\backslash B(x_0,r)$ instead of $B(x_0,r)$, where   $\Om \subset \Rn$ is  a bounded open set,  $p(x)$ is a  variable exponent, and no monotonicity type conditio is imposed onto the function  $\om(r)$ defining
 the "complementary" Morrey-type norm. In the case where $\om$ is a power function, we reveal the relation of these spaces  to weighted Lebesgue spaces. In the general case we
 prove the boundedness of the Hardy-Littlewood maximal
operator  and Calderon-Zygmund singular operators with standard kernel,  in
such spaces. We also prove  a Sobolev type ${\dual \cal M}_{\{x_0\}}^{p(\cdot),\om}
(\Om)\rightarrow {\dual \cal M}_{\{x_0\}}^{q(\cdot),\om} (\Om)$-theorem for the potential
operators $I^{\al(\cdot)},$ also of variable order.  In all the cases the
conditions  for the boundedness are given it terms of Zygmund-type integral
inequalities on $\om(r)$, which do not assume any assumption on
monotonicity of $\om(r)$.
\end{abstract}

\begin{quote}\small
\end{quote}

\begin{quote}\small
 \end{quote}

\section{Introduction}

\setcounter{theorem}{0} \setcounter{equation}{0}

In the study of local properties of solutions to of partial differential
equations, together with weighted Lebesgue spaces,  Morrey spaces
${\cal L}^{p, \lambda}(\Rn)$ play an important role, see \cite{Gi}.
Introduced by C. Morrey  \cite{Morrey} in 1938, they are defined by the norm
\begin{equation}\label{11}
\left\| f\right\|_{{\cal L}^{p,\lambda }}: = \sup_{x, \; r>0 } r^{-\frac{\lambda}{p}}
\|f\|_{L^{p}(B(x,r))} ,
\end{equation}
where $0 \le \lambda < n,$ $1\le p < \infty.$

We refer in particular to \cite{KJF} for the classical Morrey spaces. As is known, last two decades there is an increasing interest to the study of
variable exponent spaces and operators with variable parameters in such
spaces, we refer for instance to the surveying papers \cite{Din3}, \cite{K1},
\cite{KSAMADE}, \cite{Sam4} and the recent book \cite{160zc} on the progress in this field, see also references therein.

The spaces defined by the norm \eqref{11} are sometimes called \textit{global} Morrey spaces, in contrast to \textit{local} Morrey spaces
defined by the norm
\begin{equation}\label{11local}
\left\| f\right\|_{{\cal L}^{p,\lambda }_{\{x_0\}}}: = \sup_{r>0 } r^{-\frac{\lambda}{p}}
\|f\|_{L^{p}(B(x_0,r))}.
\end{equation}

 Variable exponent
Morrey spaces ${\cal L}^{p(\cdot),\lambda(\cdot)}(\Om)$, were introduced and
studied in \cite{AHS} and \cite{MizShim} in the Euclidean setting and in
\cite{KokMeskhi-Morrey} in the setting of metric measure spaces, in case of
bounded sets $\Om$. In \cite{AHS} there was proved the boundedness of the maximal
operator in variable exponent Morrey spaces ${\cal L}^{p(\cdot),\lb(\cdot)}(\Om)$  and a Sobolev-Adams type ${\cal L}^{p(\cdot),\lambda(\cdot)}\rightarrow {\cal L}^{q(\cdot),\lambda(\cdot)}$-theorem for potential operators of variable order $\al(x)$. In the case of constant $\al$, there
was also proved the ${\cal L}^{p(\cdot),\lambda(\cdot)}\to BMO$-boundedness  in the limiting case $p(x) =
\frac{n-\lambda(x)}{\al}.$ In \cite{MizShim} the
maximal operator and potential operators were considered in a somewhat more
general space, but under more restrictive conditions on $p(x)$. \ P.
H\"ast\"o in \cite{Hasto} used his "local-to-global" approach to extend
the result of \cite{AHS} on the maximal operator to the case of the whole
space $\Rn$.

  In \cite{KokMeskhi-Morrey} there was proved the boundedness of the maximal
operator and the singular integral operator in variable exponent Morrey
spaces ${\cal L}^{p(\cdot),\lambda(\cdot)}$ in the general setting of metric
measure spaces. In the case of constant $p$ and $\lb$, the results on the
boundedness of potential operators and classical Calderon-Zygmund singular
operators go back to \cite{Adams} and \cite{P}, respectively, while the
boundedness of the maximal operator in the Euclidean setting was proved in
\cite{ChiFra}; for further results  in the case of constant $p$  and $\lb$
see for instance \cite{BurGul1}-- \cite{BurGulSerTar1}.

In \cite{GulHasSam} we studied the boundedness of the classical integral operators in the generalized variable exponent  Morrey spaces \ $\M^{p(\cdot),\vi}(\Om)$ over an open bounded set $\Om \subset \Rn.$ Generalized Morrey spaces of such a kind in the case of constant $p$, with the norm  \eqref{11} replaced by
\begin{equation}\label{12}
\left\| f\right\|_{{\cal L}^{p,\vi }}: = \sup_{x, \; r>0 }\frac{r^{-\frac{n}{p}}}{\vi(r)}
\|f\|_{L^{p}(B(x,r))} ,
\end{equation}
under some assumptions on $\vi$
were
studied  in  \cite{EridGunaNakai}, \cite{KurNishSug},
\cite{Miz}, \cite{Nakai}, \cite{Nakai1}.
Results of \cite{GulHasSam}   were extended in \cite{GulHasSam1} to the case
of the generalized Morrey spaces ${\cal M}^{p(\cdot),\theta(\cdot) ,\om(\cdot)}(\Om)$ (where the $L^\infty$-norm in $r$ in the definition of the Morrey space is replaced by the Lebesgue $L^\theta$-norm), we refer to
\cite{BurHus1} for such spaces in the case of constant exponents.

 In \cite{Gul} (see, also \cite{GulBook}) there were introduced and studied local "complementary" generalized Morrey spaces ${\dual \cal M}_{\{x_0\}}^{p,\o}(\Om), \ \Om\subseteq\Rn,$ with constant  $p\in [1,\infty)$, the space of all functions $f \in L_{\textrm{loc}}^{p}(\Om\backslash \{x_0\})$ with the finite norm
\begin{equation*}\label{103}
\|f\|_{{\dual \cal M}_{\{x_0\}}^{p,\o}(\Om)}=
\sup_{r>0 } \frac{r^{\frac{n}{p^\prime}}}{\om(r)}
\|f\|_{L^{p}(\Om\backslash B(x_0,r))}.
\end{equation*}

For the particular case when $\om$ is a
power function ( \cite{Gul}, \cite{GulBook}), see also \cite{GulMus1, GulMus2}), we find it convenient to keep  the traditional notation ${\dual \cal L}_{\{x_0\}}^{p,\lambda}(\Om)$ for the space defined by the  norm
\begin{equation}\label{13}
\|f\|_{{\dual \cal L}_{\{x_0\}}^{p,\lambda}(\Om)}=
\sup_{r>0 } r^{\frac{\lambda}{p^\prime}}
\|f\|_{L^{p}(\Rn\backslash B(x_0,r))}<\infty, \ \ \ x_0\in \Om, \ \ \ 0 \le \lambda < n
\end{equation}
 Obviously, we  recover the space
${\dual \cal L}_{\{x_0\}}^{p,\lambda}(\Om)$  from ${\dual \cal M}_{\{x_0\}}^{p,\om}(\Om)$ under the choice $\om(r)=r^{\frac{n-\lambda}{p^\prime}}.$

In contrast to the Morrey space, where one measures the regularity of a function $f$ near the point $x_0$ (in the case of local Morrey spaces) and near all points $x\in\Om$ (in the case of global Morrey spaces),  the norm \eqref{13} is aimed to measure a "bad" behaviour of $f$ near the point $x_0$ in terms of the possible growth of $\|f\|_{L^{p}(\Om\backslash B(x_0,r))}$ as $r\to 0$. Correspondingly, one admits $\vi(0)=0$ in \eqref{12} and $\om(0)=\infty$ in \eqref{13}.

\vs{1mm}
  In this paper we consider local "complementary" generalized Morrey spaces ${\dual \cal M}_{\{x_0\}}^{p(\cdot),\om}(\Om)$ with variable exponent $p(\cdot)$, see Definition \ref{def0}.
However, we start with the case of constant $p$ and in this case  reveal an intimate connection of the  complementary spaces with weighted Lebesgue spaces. In the case where $\om(r)$ is a power function, we show that the space
${\dual \cal L}_{\{x_0\}}^{p,\lambda}(\Om)$
is embedded between the weighted space $L^p(\Om, |x-x_0|^{\lb(p-1)})$ and its weak version, but does not coincide with either of them, which elucidates the nature of these spaces.

 In the general case, for the spaces ${\dual \cal M}_{\{x_0\}}^{p(\cdot),\om}(\Om)$ over bounded sets $\Om \subset \Rn$ we
consider the following operators: \\
1) the Hardy-Littlewood maximal operator
$$
Mf(x)=\sup\limits_{r>0}\frac{1}{|B(x,r)|}
\intl_{\widetilde{B}(x,r)}|f(y)|dy
$$
2) variable order potential type operators
$$
I^{\al(x)} f(x)=\intl_{\Om} \frac{f(y)\,dy}{|x-y|^{n-\al(x)}}, ,
$$
3) variable order fractional maximal operator
$$
M^{\al(x)}f(x)=\sup\limits_{r>0}\frac{1}{|B(x,r)|^{1-\frac{\al(x)}{n}}}
\intl_{\widetilde{B}(x,r)}|f(y)|dy,
$$ where
$0<\inf \al(x)\le \sup\al(x)<n$,
and\\
4)  Calderon-Zygmund type singular operator
\begin{equation*}\label{s1}
T f(x)=\intl_{\Om}K(x,y) f(y)  dy
\end{equation*}
with    a "standard" singular kernel in the sense of R.Coifman and Y.Meyer, see for instance
\cite{132}, p. 99.

We find  conditions on the pair of functions  $\om_1(r)$ and $\om_2(r)$ for the  $p(\cdot)\to p(\cdot)$-boundedness of the
 operators $M$ and  $T$ from a variable exponent local "complementary" generalized
Morrey space ${\dual \cal M}_{\{x_0\}}^{p(\cdot),\om_1}(\Om)$ into another one  ${\dual \cal M}_{\{x_0\}}^{p(\cdot),\om_2}(\Om)$, and for the corresponding Sobolev $p(\cdot)\to q(\cdot)$-boundedness for the potential operators $I^{\al(\cdot)},$ under the
log-condition on $p(\cdot)$.

\vs{1mm}The paper is organized as follows.

In Section \ref{relations} we start with the case of the
spaces ${\dual \cal L}_{\{x_0\}}^{p,\lambda}(\Om)$ with constant $p$ and show their relation to the weighted Lebesgue
space $L^p(\Om, |x-x_0|^{\lb(p-1)}).$  The main statements are given in  Theorems \ref{lem78} and \ref{lem1}. In Section
\ref{Preliminaries} we
provide necessary preliminaries on variable exponent Lebesgue and Morrey
spaces. In Section \ref{generalized} we introduce the local "complementary" generalized Morrey
spaces with variable exponents and recall some facts known for generalized
Morrey spaces with constant $p$. In Section \ref{sectionmaximal} we deal with
the maximal operator, while potential  operators are studied in Section
\ref{potentials}. In Section \ref{singular} we treat Calderon-Zygmund
singular operators.

The main results are given in Theorems \ref{M1}, \ref{M1X} and   \ref{SIO1}. We
emphasize that the results we obtain for generalized Morrey spaces are new
even in the case when $p(x)$ is constant, because we do not impose any
monotonicity type condition on $\om(r).$

 \vs{6mm}
\noindent N o t a t i o n : \\
$\Rn$ is the $n$-dimensional Euclidean space,\\
 $\Om \subset \Rn$
is an open set, \ $\ell=$ diam\, $\Om$;\\
$\chi_E(x)$ is the characteristic function of a set $E\subseteq
\Rn$;\\
 $B(x,r)=\{y \in \Rn :|x-y| < r
\}), \  \widetilde{B}(x,r)=B(x,r)\cap\Om$;\\
by $c$,$C, c_1,c_2$ etc, we denote various absolute  positive constants,
which may have different values even in the same line.

\section{Relations of the "complementary" Morrey spaces ${\dual \cal L}_{\{x_0\}}^{p,\lambda}(\Om)$ with weighted Lebesgue spaces; the case of constant $p$}\label{relations}

\setcounter{theorem}{0} \setcounter{equation}{0}

 We use the standard notation $L^p(\Om, \varrho)=\left\{f: \int_\Om \varrho(y)|f(y)|^p\, dy <\infty\right\}, $
where $\varrho$ is a weight function. For the space ${\dual \cal L}_{\{x_0\}}^{p,\lambda}(\Om)$ defined in \eqref{13}, the following statement  holds.

\begin{theorem}\label{lem78} Let $\Om$ be a bounded open set, $1\le p <\infty, \ 0\le\lb\le n$ and $A>\ell$. Then
\begin{equation}\label{2a3255}
 L^p(\Om, |y-x_0|^{\lb(p-1)})    \hookrightarrow  {\dual \cal L}_{\{x_0\}}^{p,\lambda}(\Om)  \hookrightarrow
 \ \bigcap\limits_{\ve>0} L^p\left(\Om, \frac{|y-x_0|^{\lb(p-1)}}{\left(\ln\frac{A}{|y-x_0|}\right)^{1+\ve}}\right)
\end{equation}
where both the embeddings are strict, with the counterexamples $f(x)=\frac{1}{|x-x_0|^{\frac{n}{p}+\frac{\lb}{p^\prime}}}$ and $g(x)=\frac{\ln\left(\ln \frac{B}{|x-x_0|}\right)}{|x-x_0|^{\frac{n}{p}+\frac{\lb}{p^\prime}}}, \ B>\ell e^e$:
   $$f\in {\dual \cal L}_{\{x_0\}}^{p,\lambda}(\Om), \ \ \ \textrm{but} \ \ \ \ f \notin L^p(\Om, |y-x_0|^{\lb(p-1)}),$$
and
$$
g\in \bigcap\limits_{\ve>0} L^p\left(\Om, \frac{|y-x_0|^{\lb(p-1)}}{\left(\ln\frac{A}{|y-x_0|}\right)^{1+\ve}}\right), \ \ \ \ \textrm{but}
\ \ \ \ g \notin {\dual \cal L}_{\{x_0\}}^{p,\lambda}(\Om).$$
\end{theorem}

\begin{proof}\\ \textit{1$^0.\ $ The left-hand side embedding}.\\
Denote $\nu=\lb(p-1).$ For all $0<r<\ell$ we have
\begin{equation}\label{1hop}
\left(\int_\Om |y-x_0|^\nu|f(y)|^p\,dy\right)^\frac{1}{p}\ge \left(\int_{\Om\backslash \widetilde{B}(x_0,r)} |y-x_0|^\nu|f(y)|^p\,dy\right)^\frac{1}{p}\ge r^\frac{\nu}{p} \left(\int_{\Om\backslash \widetilde{B}(x_0,r)} |f(y)|^p\,dy\right)^\frac{1}{p}.
\end{equation}
Thus $\|f\|_{L^p(\Om, |y-x_0|^\nu)}\ge
 r^\frac{\lb}{p^\prime}  \|f\|_{L^p\Om\backslash \widetilde{B}(x_0,r))}$
and then
$$\|f\|_{L^p(\Om, |y-x_0|^\nu)}\ge \|f\|_{{\dual \cal L}_{\{x_0\}}^{p,\lambda}(\Om) }.$$

\textit{2$^0.\ $ The right-hand side embedding}.\\
 We take $x_0=0$ for simplicity and denote
$w_\ve(|y|)=\frac{|y|^{\lb(p-1)}}{\left(\ln\frac{A}{|y|}\right)^{1+\ve}}.$
 We have
\begin{equation}\label{trick}
\int_{\widetilde{B}(x_0,t)}|f(y)|^p w_\ve(|y|)dy  = \int_{\widetilde{B}(x_0,t)}|f(y)|^p
\left(\int_0^{|y|} \frac{d}{ds}w_\ve(s)ds\right)dy,
\end{equation}
with
$$w^\prime_\ve(t)=t^{\lb(p-1)-1}\left[\frac{\lb(p-1)}{\left(\ln\frac{A}{t}\right)^{1+\ve}}+\frac{(1+\ve)} {\left(\ln\frac{A}{t}\right)^{2+\ve}}\right]\ge 0.$$
Therefore,
$$
\int_{\widetilde{B}(x_0,t)}|f(y)|^p w_\ve(|y|)dy=  \int_0^{t} w_\ve^\prime(s) \,
\left(\int_{\{y \in \Om : s<|x_0-y|< t\}} |f(y)|^p dy\right)ds\le
$$
$$\le \int_0^\ell w^\prime_\ve(s)
\|f\|^p_{L^{p}(\Om\backslash \widetilde{B}(x_0,s))} ds.
\le \|f\|^p_{{\dual \cal M}_{\{x_0\}}^{p,\om}(\Om)}\int_0^\ell s^{\lb(p-1)}w_\ve^\prime(s)\,
 ds$$
where the last integral converges when $\ve>0$ since
$s^{\lb(p-1)}w_\ve^\prime(s)\le \frac{C}{s\left(\ln\frac{A}{s}\right)^{1+\ve}}.$ This completes the proof.
\end{proof}

\textit{3$^0.\ $ The strictness of the embeddings}.\\

Calculations for the function $f$ are obvious.  In case of the function  $g$, take $x_0=0$ for simplicity and denote   $w_\ve(|x|)=\frac{|x|^{\lb(p-1)}}{\left(\ln\frac{A}{|x|}\right)^{1+\ve}}.$ We have
$$\|f\|^p_{L^p\left(\Om,w_\ve\right)}=\int_\Om\frac{\ln^p\left(\ln \frac{B}{|x|}\right)}{|x|^n\left(\ln\frac{A}{|x|}\right)^{1+\ve}}dx \le C\intl_0^\ell
\frac{\ln^p\left(\ln \frac{B}{t}\right)}{t\left(\ln\frac{A}{t}\right)^{1+\ve}}dt<\infty$$
for every $\ve>0.$
However, for small $r\in \left(0,\frac{\dl}{2}\right)$, where $\dl=dist(0,\partial\Om)$, we obtain
$$r^\frac{\lb}{p^\prime}\|f\|_{L^p(\Om\backslash B(0,r))}= \left(r^{\lb(n-1)}\intl_{x\in\Om:\ |x|>r}
\frac{\ln^p\left(\ln \frac{B}{|x|}\right)\,dx}{|x|^{n+\lb(p-1)}}
\right)^\frac{1}{p}$$
$$\ge \left(r^{\lb(p-1)}\intl_{x\in\Om:\ r<|x| <\dl}
\frac{\ln^p\left(\ln \frac{B}{|x|}\right)\,dx}{|x|^{n+\lb(p-1)}}
\right)^\frac{1}{p}= \left(r^{\lb(p-1)}|\mathbb{S}^{n-1}|\intl_r^\dl
\frac{\ln^p\left(\ln \frac{B}{t}\right)\,dt}{t^{1+\lb(p-1)}}
\right)^\frac{1}{p}.$$
But
$$\intl_r^\dl
\frac{\ln^p\left(\ln \frac{B}{t}\right)\,dt}{t^{1+\lb(p-1)}}\ge
\intl_r^{2r}
\frac{\ln^p\left(\ln \frac{B}{t}\right)\,dt}{t^{1+\lb(p-1)}}
\ge \ln^p\left(\ln \frac{B}{2r}\right) \intl_r^{2r}
\frac{\,dt}{t^{1+\lb(p-1)}}= C \ln^p\left(\ln \frac{B}{r}\right) r^{-\lb(p-1)},
$$
so that
$$r^\frac{\lb}{p^\prime}\|f\|_{L^p(\Om\backslash B(0,r))}\ge C \ln \left(\ln \frac{B}{2r}\right) \to \infty \ \ \ \ \textrm{as} \ \ \ \ r\to 0,$$
which completes the proof of the lemma.

\begin{remark}\label{lem} The arguments similar to those  in \eqref{1hop} show that the left-hand side embedding  in \eqref{2a3255} may be extended to the case of more general
spaces ${\dual \cal M}_{\{x_0\}}^{p,\om}(\Om)$:
$$L^p(\Om, \rho(|y-x_0|))  \hookrightarrow {\dual \cal M}_{\{x_0\}}^{p,\om}(\Om)$$
where $\rho$ is a positive  increasing (or almost increasing) function such that
$\inf\limits_{r>0}\frac{\rho(r)\om^p(r)}{r^{n(p-1)}}>0. $
\end{remark}

\vs{5mm} The next theorem shows that the space  ${\dual \cal L}_{\{x_0\}}^{p,\lambda}(\Om)$ is also embedded between that weighted space $L^p(\Om, |y-x_0|^{\lb(p-1)})$ and its weak version $wL^p(\Om, |y-x_0|^{\lb(p-1)})$. The latter is defined by the norm
$$\|f\|_{wL^p(\Om, |y-x_0|^{\lb(p-1)})}=\sup_{t>0}t \left[\mu\{x\in\Om : |f(x)|>t\}\right]^\frac{1}{p}<\infty.$$
where $\mu(E)=\intl_E |y-x_0|^{\lb(p-1)}dy.$

\begin{theorem}\label{lem1}
Let $\Om$ be a bounded domain and $1\le p <\infty, \ 0<\lb\le n$. Then
\begin{equation}\label{2at7}
 L^p(\Om, |y-x_0|^{\lb(p-1)})    \hookrightarrow  {\dual \cal L}_{\{x_0\}}^{p,\lambda}(\Om)  \hookrightarrow
  wL^p(\Om, |y-x_0|^{\lb(p-1)}).
\end{equation}
\end{theorem}

\begin{proof}
By Lemma \ref{lem78}, we only have to prove  the right-hand side embedding. Let $x_0=0$. We have
$$\|f\|_{wL^p(\Om, |y|^{\lb(p-1)})}=\sup_{t>0} \left(t^p\intl_{x\in\Om : \ |f(x)|>t}|x|^{\lb(p-1)}\,dx\right)^\frac{1}{p}.$$
It suffices to estimate this norm only for small $|x|<\dl , \  \dl=dist(0,\partial\Om)$, since the embedding
$${\dual \cal L}_{\{0\}}^{p,\lambda}\left(\Om\backslash B\left(0,\dl \right)\right)  \hookrightarrow
  L^p\left(\Om\backslash B\left(0,\dl \right), |y|^{\lb(p-1)}\right) \hookrightarrow
  wL^p\left(\Om\backslash B\left(0,\dl \right), |y|^{\lb(p-1)}\right)$$
  is obvious.
  We obtain
$$\sup_{t>0}t^p \intl_{x\in B\left(0,\dl \right) : \ |f(x)|>t}|x|^{\lb(p-1)}\,dx =
\sup_{t>0}t^p \intl_{x\in B\left(0,\dl \right) : \ |f(x)|>t}\frac{\left\{|x|^\frac{\lb}{p^\prime}
\|f\|_{L^p(\Om\backslash B(0,|x|))}
\right\}^p}{\intl_{\Om\backslash B(0,|x|)}|f(y)|^p\,dy}\,dx$$
Hence
$$\|f\|_{wL^p(\Om, |y|^{\lb(p-1)})}\le \|f\|_{{\dual \cal L}_{\{x_0\}}^{p,\lambda}(\Om)}\left(
\intl_{B(0,\dl)}\frac{dx}{|\Om\backslash B(0,|x|)|}\right) ^\frac{1}{p}
\le \left(
\frac{|B(0,\dl)|}{|\Om\backslash B(0,\dl)|}\right) ^\frac{1}{p}\|f\|_{{\dual \cal L}_{\{x_0\}}^{p,\lambda}(\Om)}
$$
which proves the lemma.
\end{proof}

\section{Preliminaries on variable exponent Lebesgue and Morrey
spaces}\label{Preliminaries}

\setcounter{theorem}{0} \setcounter{equation}{0}

We refer to the book \cite{160zc} for variable exponent Lebesgue spaces but give some basic definitions and facts. Let $p(\cdot)$ be a measurable function on $\Om$ with values
in $[1,\infty)$. An open set $\Om$ is assumed to be bounded
throughout the whole paper. We mainly suppose that
\begin{equation}\label{h0}
1< p_-\le p(x)\le p_+<\infty,
\end{equation}
where $
 p_-:=\ei_{x\in \Om}p(x),  \ p_+:=\es_{x\in
\Om}p(x).
$
By $L^{p(\cdot)}(\Om)$ we denote the space of all measurable functions $f(x)$
on $\Om$ such that
$$
I_{p(\cdot)}(f)= \int_{\Om}|f(x)|^{p(x)} dx < \infty.
$$
Equipped with the norm
$$\|f\|_{p(\cdot)}=\inf\left\{\eta>0:~I_{p(\cdot)}\left(\frac{f}{\eta}\right)\le
1\right\},
$$
this is a Banach function space. By $p^\prime(x) =\frac{p(x)}{p(x)-1}$, $x\in
\Om,$ we denote the conjugate exponent.

By $WL(\Om)$ (weak Lipschitz) we denote the class of functions defined on
$\Om$ satisfying the log-condition
\begin{equation}\label{h8}
|p(x)-p(y)|\le\frac{A}{-\ln |x-y|}, \;\; |x-y|\le \frac{1}{2} \;\; x,y\in \Om,
\end{equation}
where  $A=A(p) > 0$ does not depend on $x, y$.

\begin{theorem}\label{D} (\cite{Din}) Let $\Om \subset \Rn$ be an open bounded set  and $p\in WL(\Om)$ satisfy
condition \eqref{h0}. Then the maximal operator $M$ is bounded in
$L^{p(\cdot)}(\Om)$.
\end{theorem}

\vs{5mm} The following theorem  was proved in \cite{Sam2} under the condition
that the maximal operator is bounded in $L^{p(\cdot)}(\Om)$, which
became an unconditional result after the result of Theorem
\ref{D}.
\begin{theorem}\label{S1} Let $\Om \subset \Rn$ be bounded, $p,\al\in WL(\Om)$ satisfy
assumption \eqref{h0} and the conditions
\begin{equation}\label{1}
\inf\limits_{x\in \Om} \al(x)>0,\,\, \sup\limits_{x\in \Om} \al(x)p(x)<n.
\end{equation}
 Then the operator $I^{\al(\cdot)}$ is bounded from $L^{p(\cdot)}(\Om)$ to
$L^{q(\cdot)}(\Om)$ with $
\frac
1{q(x)}=\frac 1{p(x)}-\frac {\al(x)} {n}.$
\end{theorem}

\vs{2mm}
Singular operators within the framework of the spaces with variable exponents
were studied in \cite{Din2}. From Theorem 4.8 and Remark 4.6 of \cite{Din2}
and the known results on the boundedness of the maximal operator, we have the
following statement, which is
  formulated below for our goals for a bounded $\Om$, but valid for an arbitrary open set $\Om$
under the corresponding condition on $p(x)$ at infinity.

\begin{theorem}\ (\cite{Din2})\label{SIO} Let $\Om \subset \Rn$ be a bounded
 open set and  $p\in WL(\Om)$ satisfy
condition \eqref{h0}. Then the singular integral operator $T$ is
bounded in $L^{p(\cdot)}(\Om)$.
\end{theorem}

The estimate provided by the following lemma (see
\cite{Sam2}, Corollary to Lemma 3.22) is crucial for our further proofs.

\begin{lemma}\label{lemma} Let $\Om$ be a bounded domain and $p$
satisfy the  condition
\eqref{h8} and  $1\le p_-\le p(x)\le p_+<\infty$. Let also $\sup_{x\in\Om} \nu(x)<\infty$ and $\sup_{x\in\Om}[n+\nu(x)p(x)]<0$. Then
\begin{equation}\label{estikmate}
\||x-\cdot|^{\nu(x)}\chi_{\Om\backslash \widetilde{B}(x,r)}(\cdot)\|_{p(\cdot)}\le C
r^{\nu(x)+\frac{n}{p(x)}}, \quad \ x\in\Om, 0<r<\ell=diam\,\Om,
\end{equation}
where $C$ does not depend on $x$ and $r$.
\end{lemma}

\vs{4mm} Let $ \lambda(x)$ be a measurable function on $\Om$ with values in
$[0,n]$. The variable Morrey space ${\cal L}^{p(\cdot),\lambda(\cdot)}(\Om)$ introduced in \cite{AHS},
is defined   as the set of integrable functions $f$ on $\Om$ with the finite
norm
$$
\|f\|_{{\cal L}^{p(\cdot),\lambda(\cdot)}(\Om)}= \sup_{x\in \Om, \; t>0}
t^{-\frac{\lambda(x)}{p(x)}}\|
f\chi_{\widetilde{B}(x,t)}\|_{L^{p(\cdot)}(\Om)}.
$$

The following statements are known.

\begin{theorem}\label{maximal} (\cite{AHS}) Let $\Om$ be bounded and $p\in WL(\Om)$ satisfy
condition \eqref{h0} and let a measurable function $\lb$ satisfy
the conditions
$
0\le \lb(x), \ \quad \sup_{x\in \Om} \lambda(x)<n.
$
 Then the maximal operator $M$ is bounded in ${\cal L}^{p(\cdot),\lambda(\cdot)} (\Om)$.
\end{theorem}

Theorem \ref{maximal} was extended to unbounded domains in \cite{Hasto}.  Note that the boundedness of the maximal operator
  in Morrey spaces with variable $p(x)$ was
studied  in \cite{KokMeskhi-Morrey}  in the more general setting of
quasimetric measure spaces.

The known statements for the potential operators read as follows.

\begin{theorem}\label{Spanne-Ris} (\cite{AHS}; Spanne-type result).  Let $\Om$ be bounded, $p,\al
\in WL(\Om)$ and $p$ satisfy condition \eqref{h0}. Let also $\lb(x)\ge 0$,
the conditions \eqref{1} be fulfilled and  $\frac
1{q(x)}=\frac 1{p(x)}-\frac {\al(x)} {n}.$
Then the operator $I^{\al(\cdot)}$ is bounded from ${\cal L}^{p(\cdot),\lambda(\cdot)} (\Om)$ to ${\cal L}^{q(\cdot),\mu(\cdot)}
(\Om)$, where
$
\frac{\mu(x)}{q(x)}=\frac{\lambda(x)}{p(x)}.
$
\end{theorem}

\begin{theorem}\label{Ris} (\cite{AHS}; Adams-type result). Let $\Om$ be bounded, $p,\al
\in WL(\Om)$ and $p$ satisfy condition \eqref{h0}. Let also $\lb(x)\ge 0$ and
the conditions
\begin{equation}\label{2}
\inf\limits_{x\in \Om} \al(x)>0,\,\, \sup\limits_{x\in \Om}
[\lambda(x)+\al(x)p(x)]<n
\end{equation}
hold. Then the operator $I^{\al(\cdot)}$ is bounded from ${\cal L}^{p(\cdot),\lambda(\cdot)} (\Om)$ to
${\cal L}^{q(\cdot),\lambda(\cdot)}
(\Om)$, where
$ \frac 1{q(x)}=\frac 1{p(x)}-\frac {\al} {n-\lambda(x)}.
$
\end{theorem}

\section{Variable exponent local "complementary" generalized Morrey spaces}\label{generalized}
\setcounter{theorem}{0} \setcounter{equation}{0}

Everywhere in the sequel the functions $\om(r),\ \om_1(r)$ and
$\om_2(r)$ used in the body of the paper, are non-negative measurable
function on $(0,\ell)$, $\ell=diam\,\Om$. Without loss of generality we may assume that they are bounded beyond any small neighbourhood $(0,\dl)$ of the origin.

The local generalized Morrey space
${\cal M}_{\{x_0\}}^{p(\cdot),\om}(\Om)$ and global generalized Morrey spaces
${\cal M}^{p(\cdot),\om}(\Om)$  with variable exponent are  defined  (see \cite{GulHasSam}) by the norms
$$
\|f\|_{{\cal M}_{\{x_0\}}^{p(\cdot),\om}} =
\sup\limits_{r>0}\frac{r^{-\frac{n}{p(x_0)}}}{\om(r)}
\|f\|_{L^{p(\cdot)}(\widetilde{B}(x_0,r))},
$$
$$
\|f\|_{{\cal M}^{p(\cdot),\om}} =
\sup\limits_{x \in \Om, r>0}\frac{r^{-\frac{n}{p(x)}}}{\om(r)}
\|f\|_{L^{p(\cdot)}(\widetilde{B}(x,r))},
$$
where  $x_0 \in \Om$ and $1 \le p(x) \le p_+< \infty$ for all $x \in \Om.$

\vs{2mm}
We find it convenient to introduce the variable exponent version of the local  "complementary" space as follows (compare with  \eqref{13}).
\begin{definition}\label{def0}
 Let  $x_0 \in \Om$, $1 \le p(x) \le p_+< \infty$.
The  \textit{variable exponent generalized local "complementary"  Morrey space}
${\dual \cal M}_{\{x_0\}}^{p(\cdot),\om}(\Om)$
 is defined by the norm
$$
\|f\|_{{\dual \cal M}_{\{x_0\}}^{p(\cdot),\om}} =
\sup\limits_{r>0}\frac{r^{\frac{n}{p^\prime(x_0)}}}{\om(r)}
\|f\|_{L^{p(\cdot)}(\Om\backslash \widetilde{B}(x_0,r))}.
$$
\end{definition}

Similarly to the notation in \eqref{13}, we introduce the following particular case of the space  ${\dual \cal M}_{\{x_0\}}^{p(\cdot),\om}(\Om),$ defined by the norm
\begin{equation}\label{13suo}
\|f\|_{{\dual \cal L}_{\{x_0\}}^{p(\cdot),\lambda}(\Om)}=
\sup_{r>0 } r^{\frac{\lambda}{p^\prime}}
\|f\|_{L^{p(\cdot)}(\Om\backslash B(x_0,r))}<\infty, \ \ \ x_0\in \Om, \ \ \ 0 \le \lambda < n
\end{equation}

Everywhere in the sequel we assume that
\begin{equation} \label{BHF}
\sup\limits_{0<r<\ell}\frac{r^{\frac{n}{p^\prime(x_0)}}}{\om(r)}<\infty , \ \ \ \ \ell=\diam\ \Om,
 \end{equation}
which makes the space ${\dual \cal M}_{\{x_0\}}^{p(\cdot),\om}(\Om)$ non-trivial, since it contains $L^{p(\cdot)}(\Om)$ in this case.

\begin{remark}\label{rem1} Suppose that
$$\inf\limits_{\dl<r<\ell}\frac{r^{\frac{n}{p^\prime(x_0)}}}{\om(r)}>0$$
for every $\dl>0.$ Then
$$\|f\|_{{\dual \cal M}_{\{x_0\}}^{p(\cdot),\om}(\Om)}\sim \|f\|_{{\dual \cal M}_{\{x_0\}}^{p(\cdot),\om}(B(x_0,\dl))}+ \|f\|_{L^{p(\cdot)}(\Om\backslash B(x_0,\dl))}$$
(with the constants in the above equivalence depending on $\dl$). Since $\dl>0$ is arbitrarily small, the space
 ${\dual \cal M}_{\{x_0\}}^{p(\cdot),\om}(\Om)$
 may be  interpreted as the space of functions    \textit{whose restrictions onto an arbitrarily small neighbourhood $B(x_0,\dl)$  are in local "complementary" Morrey space ${\dual \cal M}_{\{x_0\}}^{p(\cdot),\om}(B(x_0,\dl))$ with the  exponent $p(\cdot)$ close to the constant value $p(x_0)$ and the restrictions onto the exterior $\Om\backslash B(x_0,\dl)$ are
 in the  variable exponent Lebesgue space $L^{p(\cdot)}$.}

 \end{remark}

\vs{2mm}
 If also $\inf\limits_{0<r<\ell}\frac{r^{\frac{n}{p^\prime(x_0)}}}{\om(r)}>0$, then
${\dual \cal M}_{\{x_0\}}^{p(\cdot),\om}(\Om)=L^{p(\cdot)}(\Om)$. Therefore, to guarantee that the "complementary" space ${\dual \cal M}_{\{x_0\}}^{p(\cdot),\om}(\Om)$ is strictly larger than $L^{p(\cdot)}(\Om),$ one should be interested in the cases where
 \begin{equation} \label{14}
\lim\limits_{r\to 0}\frac{r^{\frac{n}{p^\prime(x_0)}}}{\om(r)}=0.
 \end{equation}

 Clearly,  the space ${\dual \cal M}_{\{x_0\}}^{p(\cdot),\om}(\Om)$ may contain functions with a non-integrable singularity at the point $x_0$, if no additional assumptions are introduced.
To study the operators in ${\dual \cal M}_{\{x_0\}}^{p(\cdot),\om}(\Om),$
we need its embedding into $L^1(\Om).$  Corollary below shows that the Dini condition on $\om$
is sufficient for such an embedding.

\begin{lemma}\label{integrability}
Let $f\in L^{p(\cdot)}(\Om\backslash \widetilde{B}(x_0,s))$ for every $s\in (0.\ell)$ and $\gm\in \mathbb{R}.$ Then
\begin{equation} \label{Ga14}
\int_{\widetilde{B}(x_0,t)}|y-x_0|^\gm |f(y)| dy \le C \int_0^{t} s^{\gm+\frac{n}{p^\prime(x_0)} -1}
\|f\|_{L^{p(\cdot)}(\Om\backslash \widetilde{B}(x_0,s))}ds
\end{equation}
for every $t\in (0,\ell),$ where $C$ does not depend on $f,t$ and $x_0$.
\end{lemma}
\begin{proof}
We use the following trick, where the
 parameter $\bt>0$ which will be chosen  later:
\begin{equation}\label{trick}
\int_{\widetilde{B}(x_0,t)}|y-x_0|^\gm|f(y)| dy  = \beta\int_{\widetilde{B}(x_0,t)}|x_0-y|^{\gm-\beta}|f(y)|
\left(\int_0^{|x_0-y|} s^{\beta-1}ds\right)dy
\end{equation}
$$
=\beta \int_0^{t} s^{\beta-1} \,
\left(\int_{\{y \in \Om : s<|x_0-y|< t\}}|x_0-y|^{\gm -\beta} \, |f(y)| dy\right)ds.
$$
Hence
$$\int_{\widetilde{B}(x_0,t)}|f(y)| dy\le C \int_0^{t} s^{\beta-1} \,
\|f\|_{L^{p(\cdot)}(\Om\backslash \widetilde{B}(x_0,s))} \; \||x_0-y|^{\gm -\beta}\|_{L^{p^\prime(\cdot)}(\Om\backslash \widetilde{B}(x_0,s))} ds.
$$
Now we make use of Lemma  \ref{lemma} which is possible if we choose $\beta >\max\left(0,\frac {n}{p^\prime_-}+\gm\right)$ and then arrive at
 \eqref{Ga14}.
 \end{proof}

\begin{corollary}\label{cor}
The following embeddings hold
\begin{equation}\label{embeddungs}
L^{p(\cdot)}(\Om)\hookrightarrow {\dual \cal M}_{\{x_0\}}^{p(\cdot),\om}(\Om) \hookrightarrow L^1(\Om)
\end{equation}
where the  left-hand side embedding is guaranteed by the condition \eqref{BHF} and the right-hand side one by the condition
\begin{equation}\label{convergence}
\intl_0^\ell \frac{\om(r)\, dr}{r}<\infty.
\end{equation}
\end{corollary}
\begin{proof}
The statement for the left-hand side embedding is obvious. The right-hand side follows from Lemma \ref{integrability} with $\gm=0$ and
 the inequality
\begin{equation}\label{qwet}
\int_{0}^\ell r^{\frac{n}{p^\prime(x_0)}-1} \|f\|_{L^{p(\cdot)}(\Om\backslash \widetilde{B}(x_0,r))} dr
\le \|f\|_{{\dual \cal M}_{\{x_0\}}^{p(\cdot),\om}(\Om)} \int_0^\ell \om(r)\frac{dr}{r}.
\end{equation}
\end{proof}

\begin{remark}\label{rem1a}
Note that similarly to the arguments in the proof of Corollary \ref{cor}, we can see that the condition
$
\intl_0^\ell \om(r)r^{\gm-1}\, dr<\infty,  \ \gm\in \mathbb{R},
$
guarantees the embedding
$
{\dual \cal M}_{\{x_0\}}^{p(\cdot),\om}(\Om) \hookrightarrow L^1(\Om,|y-x_0|^\gm)$ into the weighted space; note that  only the values $\gm>-n/p^\prime(x_0)$ may be of interest for us, because the above condition
with
$\gm\le-n/p^\prime(x_0)$   is not compatible with the condition \eqref{BHF} of the non-triviality of the  space ${\dual \cal M}_{\{x_0\}}^{p(\cdot),\om}(\Om).$
  \end{remark}

\begin{remark}\label{rem2}
We also find it convenient to give a condition for  $f\in L^1(\Om)$ in the form
as follows:
\begin{equation}\label{dobavleno}
\int_0^{\ell} t^{\frac{n}{p^\prime(x_0)}-1}
\|f\|_{L^{p(\cdot)}(\Om\backslash \widetilde{B}(x_0,t))}dt<\infty \  \ \ \ \Longrightarrow \ \ \ \ f\in L^1(\Om),
\end{equation}
for which it suffices to refer to \eqref{Ga14}.
\end{remark}

In the sequel all the operators under consideration (maximal, singular and potential operators) will be considered on function $f$ either satisfying the condition of the existence of the integral in \eqref{dobavleno}, or belonging to
${\dual \cal M}_{\{x_0\}}^{p(\cdot),\om}(\Om)$ with $\om$ satisfying the condition \eqref{convergence}. Such functions are therefore integrable on $\Om$
in both the cases, and consequently all the studied operators exist on such functions
a.e.

\vs{3mm}Note that the  statements on the boundedness  of  the maximal, singular and potential operators in the "complementary" Morrey spaces known for the case of the constant exponent $p$, obtained in \cite{Gul}, read as follows. Note that the theorems below do not assume no monotonicity type conditions on the functions $\om, \om_1$ and $\om_2$.

\begin{theorem} (\cite{Gul}, Theorem 1.4.6) \label{nakaiVagif1}
Let $1 < p < \infty$, $x_0 \in \Rn$ and  $\om_1(r)$ and $\om_2(r)$ be positive
measurable functions satisfying the condition
 \begin{equation*} 
\int_0^r \om_1(t)\frac{dt}{t} \le  c \,\om_2(r)
 \end{equation*}
with  $c>0 $  not  depending on $r>0$. Then the
operators $M$ and $T$ are bounded from
${\dual \cal M}_{\{x_0\}}^{p,\om_1}(\Rn)$ to
${\dual \cal M}_{\{x_0\}}^{p,\om_2}(\Rn)$.
\end{theorem}
\begin{corollary} \label{Vagif1} (\cite{Gul})
Let $1 < p < \infty$, $x_0 \in \Rn$ and  $0 \le \lambda < n$.
Then the operators $M$ and $T$ are bounded in the space
${\dual \cal L}_{\{x_0\}}^{p,\lambda}(\Rn)$.
\end{corollary}

\begin{theorem}\label{Guliev} (\cite{Gul}, Theorem 1.3.9) \label{nakaiVagif2}
Let $0<\a<n$, $1 < p < \infty$, $\frac{1}{q}=\frac{1}{p}-\frac{\a}{n}$,
$x_0 \in \mathbb{R}^n$ and  $\om_1(r)$, $\om_2(r)$ be
positive measurable functions satisfying the condition
 \begin{equation*}
r^\a  \int_0^r \om_1(t)\frac{dt}{t} \le  c \,\om_2(r),
 \end{equation*}
with  $c>0 $  not  depending on $r>0$. Then the operators $M^{\a}$ and $I^{\a}$ are bounded from
${\dual \cal M}_{\{x_0\}}^{p,\om_1}(\Rn)$ to
${\dual \cal M}_{\{x_0\}}^{q,\om_2}(\Rn)$.
\end{theorem}
\begin{corollary}  \label{Vagif2}(\cite{Gul}) Let $0<\a<n$, $1<p<\frac{n}{\a}$, $x_0 \in \mathbb{R}^n$ and
$\frac{1}{p}-\frac{1}{q}=\frac{\a}{n}$ and $\frac{\lambda}{p^\prime}=\frac{\mu}{q^\prime}$.
Then the operators $M^{\a}$ and $I^{\a}$ are bounded from
${\dual \cal L}_{\{x_0\}}^{p,\lambda}(\Rn)$ to
${\dual \cal L}_{\{x_0\}}^{q,\mu}(\Rn)$.
\end{corollary}

\begin{remark}\label{rem}
The introduction of global "complementary"  Morrey-type spaces has no big sense, neither in case of constant exponents, nor in the case of variable exponents. In the case of constant exponents this was noted in  \cite{69acb}, pp 19-20; in this case the global space defined by the norm
$$\sup\limits_{x\in\Om, r>0}\frac{r^{\frac{n}{p}}}{\om(r)}
\|f\|_{L^{p}(\Om\backslash \widetilde{B}(x,r))}$$
reduces to $L^p(\Om)$ under the assumption \eqref{BHF}. In the case of variable exponents there happens the same. In general, to make it clear, note that for instance under the assumption \eqref{14} if we admit that  $\sup\limits_{ r>0}\frac{r^{\frac{n}{p(\cdot)}}}{\om(r)}
\|f\|_{L^{p}(\Om\backslash \widetilde{B}(x,r))}$ for two different points $x=x_0$ and $x=x_1, \ x_0\ne x_1$, this would immediately imply that $f\in L^{p(\cdot)}$ in a neighbourhood of both the points $x_0$ and $x_1$.
 \end{remark}

\section{The maximal operator in the spaces ${\dual \cal M}_{\{x_0\}}^{p(\cdot),\om}
(\Om)$}\label{sectionmaximal} \setcounter{theorem}{0}
\setcounter{equation}{0}

The proof of the main result of this section presented in Theorem  \ref{M1} is based on the  estimate given in the following preliminary theorem.

\begin{theorem} \label{HG21} Let $\Om$ be bounded, $p\in WL(\Om)$ satisfy the
condition \eqref{h0} and  $f \in L^{p(\cdot)}(\Om\backslash \widetilde{B}(x_0,r))$ for every $r\in (0,\ell)$. If the integral
\begin{equation}\label{converg}
\int_{0}^{\ell}
r^{\frac{n}{p^\prime(x_0)}-1}\|f\|_{L^{p(\cdot)}(\Om\backslash \widetilde{B}(x_0,r))} dr
\end{equation}
converges, then
\begin{equation}\label{GOP}
\|Mf\|_{L^{p(\cdot)}(\Om\backslash \widetilde{B}(x_0,t))} \le
Ct ^{-\frac{n}{p^\prime(x_0)}}\int_{0}^{t}
r^{\frac{n}{p^\prime(x_0)}-1}\|f\|_{L^{p(\cdot)}(\Om\backslash \widetilde{B}(x_0,r))} dr
\end{equation}
for every $t\in (0,\ell),$ where $C$ does not depend on
 $f, t$ and $x_0$.
\end{theorem}

\begin{proof}
We represent  $f$ as
\begin{equation}\label{repr}
f=f_1+f_2, \ \quad f_1(y)=f(y)\chi _{\Om\backslash \widetilde{B}(x_0,t)}(y) \quad f_2(y)=f(y)\chi _{\widetilde{B}(x_0,t)}(y).
\end{equation}
\textit{$1^o.$ Estimation of $Mf_1$.} This case is easier, being treated by means of  Theorem \ref{D}.
Obviously $f_1\in L^{p(\cdot)}(\Om)$ so that
by  Theorem \ref{D}
\begin{equation} \label{kkk}
\|Mf_1\|_{L^{p(\cdot)}(\Om\backslash \widetilde{B}(x_0,t))} \le
\|Mf_1\|_{L^{p(\cdot)}(\Om)}
\le C \|f_1\|_{L^{p(\cdot)}(\Om)} = C\|f\|_{L^{p(\cdot)}(\Om\backslash \widetilde{B}(x_0,t))}.
\end{equation}

By the monotonicity of the norm $\|f\|_{L^{p(\cdot)}(\Om\backslash \widetilde{B}(x_0,r))}$  with respect to $r$ we have
\begin{equation}\label{ADI}
\|f\|_{L^{p(\cdot)}(\Om\backslash \widetilde{B}(x_0,t))}\le Ct^{-\frac{n}{p^\prime(x_0)}} \int_{0}^{t}
r^{\frac{n}{p^\prime(x_0)}-1}\|f\|_{L^{p(\cdot)}(\Om\backslash \widetilde{B}(x_0,r))} dr
\end{equation}
and then
\begin{equation}\label{Ga12'}
\|Mf_1\|_{L^{p(\cdot)}(\Om\backslash \widetilde{B}(x_0,t))}    \le Ct^{-\frac{n}{p^\prime(x_0)}}
\int_{0}^{t} r^{\frac{n}{p^\prime(x_0)}-1}\|f\|_{L^{p(\cdot)}(\Om\backslash \widetilde{B}(x_0,r))} dr. \notag
\end{equation}

\textit{$2^o.$ Estimation of $Mf_2$.} This case needs the application of
Lemma \ref{integrability}.
To estimate $Mf_2(z)$ by means of \eqref{Ga14}, we  observe that for $z\in \Om\backslash \widetilde{B}(x_0,2t)$  we have
\begin{align*}
Mf_2(z) &= \sup\limits_{r>0}|B(z,r)|^{-1}\int_{\widetilde{B}(z,r)}|f_2(y)|dy
\\
& \le  \sup\limits_{r \ge t} \int_{ \widetilde{B}(x_0,t)\cap B(z,r)}
|y-z|^{-n} \, |f(y)|dy
\\
& \le 2^{n} \, |x_0-z|^{-n} \, \int_{ \widetilde{B}(x_0,t)}  |f(y)|dy.
\end{align*}

Then by \eqref{Ga14}
\begin{equation}\label{1001}
Mf_2(z) \le C |x_0-z|^{-n} \, \int_0^{t} s^{\frac{n}{p^\prime(x_0)}-1}
\|f\|_{L^{p(\cdot)}(\Om\backslash \widetilde{B}(x_0,s))}ds.
\end{equation}
Therefore
\begin{align}\label{Ga13'}
\|Mf_2\|_{L^{p(\cdot)}(\Om\backslash \widetilde{B}(x_0,2t))} &\le
C \int_0^{t} s^{\frac{n}{p^\prime(x_0)}-1}
\|f\|_{L^{p(\cdot)}(\Om\backslash \widetilde{B}(x_0,s))}ds\,\left\||x_0-z|^{-n}\right\|_{L^{p(\cdot)}(\Om\backslash \widetilde{B}(x_0,2t))}\notag
\\
&\le C t^{-\frac{n}{p^\prime(x_0)}}\int_0^{t} s^{\frac{n}{p^\prime(x_0)}-1}
\|f\|_{L^{p(\cdot)}(\Om\backslash \widetilde{B}(x_0,s))}ds.
\end{align}

Since $
\|Mf\|_{L^{p(\cdot)}(\Om\backslash \widetilde{B}(x_0,2t))}
\le \|Mf_1\|_{L^{p(\cdot)}(\Om\backslash \widetilde{B}(x_0,2t))}
+\|Mf_2\|_{L^{p(\cdot)}(\Om\backslash \widetilde{B}(x_0,2t))},
$ from \eqref{Ga12'} and \eqref{Ga13'} we arrive at \eqref{GOP} with $\|Mf\|_{L^{p(\cdot)}(\Om\backslash \widetilde{B}(x_0,2t))}$ on the left-hand side and then \eqref{GOP} obviously holds also for $\|Mf\|_{L^{p(\cdot)}(\Om\backslash \widetilde{B}(x_0,t))}$.
\end{proof}

 The following theorem for the complementary Morrey spaces is, in a sense,  a  counterpart to Theorem \ref{maximal} formulated in Section \ref{Preliminaries} for the usual Morrey spaces.

\begin{theorem}\label{M1}
Let $\Om \subset \Rn$ be an open bounded  set,  $p\in WL(\Om)$ satisfy the
assumption \eqref{h0} and the functions $\om_1(t)$ and $\om_2(t)$
satisfy the condition
\begin{equation}\label{eq3.6.VZ}
\int_0^{t} \om_1(r)\frac{dr}{r} \le
 C \,\om_2(t),
\end{equation}
where $C$ does not depend on $t$. Then the maximal operator $M$ is
bounded from the space ${\dual \cal M}_{\{x_0\}}^{p(\cdot),\om_1} (\Om)$ to the space
${\dual \cal M}_{\{x_0\}}^{p(\cdot),\om_2}(\Om)$.
\end{theorem}

\begin{proof} For $f\in {\dual \cal M}_{\{x_0\}}^{p(\cdot),\om_1}(\Om)$ we have
$$
\|Mf\|_{{\dual \cal M}_{\{x_0\}}^{p(\cdot),\om_2}(\Om)}=\sup_{t\in (0,\ell)}
\frac{t^\frac{n}{p^\prime(x_0)}}{\om_2(t)}
\|Mf\|_{L^{p(\cdot)}(\Om\backslash \widetilde{B}(x_0,t))},
$$
where  Theorem \ref{HG21} is applicable to the norm $\|Mf\|_{L^{p(\cdot)}(\Om\backslash \widetilde{B}(x_0,t))}$.
Indeed from \eqref{eq3.6.VZ} it follows that
the integral $\int_0^t\frac{\om_1(r)}{r}dr$ converges. This implies that  for $f\in {\dual \cal M}_{\{x_0\}}^{p(\cdot),\om_1} (\Om)$
the assumption of the convergence of the integral of type \eqref{converg} is fulfilled by \eqref{qwet}.
Then by Theorem \ref{HG21} we obtain
\begin{align*}
\|Mf\|_{{\dual \cal M}_{\{x_0\}}^{p(\cdot),\om_2}(\Om)} & \le C
\sup_{0<t\le \ell}\om^{-1}_2(t)\int_0^{t}
r^{-\frac{n}{p^\prime(x_0)}-1}\|f\|_{L^{p(\cdot)}(\Om\backslash \widetilde{B}(x_0,r))} dr.
\end{align*}
Hence
$$ \|Mf\|_{{\dual \cal M}_{\{x_0\}}^{p(\cdot),\om_2}(\Om)}  \le
C \|f\|_{{\dual \cal M}_{\{x_0\}}^{p(\cdot),\om_1}(\Om)}
\sup_{t\in(0,\ell)}\frac1{\om_2(t)} \int_0^{t} \om_1(r)\frac{dr}{r}
\le C\|f\|_{{\dual \cal M}_{\{x_0\}}^{p(\cdot),\om_1}(\Om)}
$$
by \eqref{eq3.6.VZ}, which completes the proof.
\end{proof}

\begin{corollary} \label{VHS1}
 Let $\Om \subset \Rn$ be an open bounded  set,  $x_0 \in \Om$, $0 \le \lambda < n$, $\lambda \le \mu \le n$ and let  $p\in WL(\Om)$ satisfy the
assumption \eqref{h0}.
Then the operator $M$ is bounded from the space ${\dual \cal L}_{\{x_0\}}^{p(\cdot),\lambda}(\Om)$
to ${\dual \cal L}_{\{x_0\}}^{p(\cdot),\mu}(\Om)$.
\end{corollary}

\

\section{Riesz potential operator in the spaces
 ${\dual \cal M}_{\{x_0\}}^{p(\cdot),\om}(\Om)$}\label{potentials}
\setcounter{theorem}{0} \setcounter{equation}{0}

\

In this section we  extend Theorem \ref{Guliev} to the variable exponent
setting.
  Note that Theorems \ref{G4} and
\ref{M1X} in the case of constant exponent $p$   were proved in
\cite{Gul},  Theorems 1.3.2 and 1.3.9 (see also \cite{GulBook}, p. 112, 129).

\begin{theorem} \label{G4} Let \eqref{h0} be fulfilled and   $p(\cdot), \al(\cdot) \in WL(\Om)$   satisfy the conditions in \eqref{1}.
If $f$ is such that the integral \eqref{converg}
converges, then
\begin{equation}\label{GOPRiesz}
\|I^{\al(\cdot)} f\|_{L^{q(\cdot)}(\Om\backslash \widetilde{B}(x_0,t))} \le
Ct^{-\frac{n}{p^\prime(x_0)}}\int_0^{t} s^{\frac{n}{p^\prime(x_0)}-1}\|f\|_{L^{p(\cdot)}(\Om\backslash \widetilde{B}(x_0,s))} ds
\end{equation}
for every $f \in L^{p(\cdot)}(\Om\backslash \widetilde{B}(x_0,t))$, where
\begin{equation}\label{1cxz}
\frac{1}{q(x)}=\frac 1{p(x)}-\frac {\al(x)} {n}.
\end{equation} and $C$ does not depend on
 $f, x_0$ and $ t\in (0,\ell)$.
\end{theorem}

\begin{proof}
We  represent the function $f$ in the form $f=f_1+f_2$ as in \eqref{repr} so that

$$
\|I^{\a(\cdot)} f\|_{L^{q(\cdot)}(\Om\backslash \widetilde{B}(x_0,2t))} \le
\|I^{\a(\cdot)} f_1\|_{L^{q(\cdot)}(\Om\backslash \widetilde{B}(x_0,2t))}+
\|I^{\a(\cdot)} f_2\|_{L^{q(\cdot)}(\Om\backslash \widetilde{B}(x_0,2t))}.
$$
Since $f_1\in L^{p(\cdot)}(\Om)$,  by  Theorem \ref{S1} we have
\begin{gather*}
\|I^{\a(\cdot)} f_1\|_{L^{q(\cdot)}(\Om\backslash \widetilde{B}(x_0,2t))} \le
\|I^{\a(\cdot)} f_1\|_{L^{q(\cdot)}(\Om)}
\le C \|f_1\|_{L^{p(\cdot)}(\Om)}=C \|f\|_{L^{p(\cdot)}(\Om\backslash \widetilde{B}(x_0,t))}
\end{gather*}
and then
\begin{equation}\label{Ha1}
\|I^{\a(\cdot)} f_1\|_{L^{q(\cdot)}(\Om\backslash \widetilde{B}(x_0,2t))} \le C t^{-\frac{n}{p^\prime(x_0)}} \int_0^{t} s^{\frac{n}{p^\prime(x_0)}-1}\|f\|_{L^{p(\cdot)}(\Om\backslash \widetilde{B}(x_0,s))}ds
\end{equation}
in view of \eqref{ADI}.

To estimate
$$\|I^{\a(\cdot)} f_2\|_{L^{q(\cdot)}(\Om\backslash \widetilde{B}(x_0,2t))}=
\left\|\int\limits_{\widetilde{B}(x_0,t)}|z-y|^{\al(z) -n}f(y) dy\right\|_{L^{q(\cdot)}(\Om\backslash \widetilde{B}(x_0,2t))},
$$
we observe that for $z\in \Om\backslash \widetilde{B}(x_0,2t)$ and $y\in \widetilde{B}(x_0,t)$ we have $\frac{1}{2} |x_0-z| \le |z-y|\le\frac{3}{2} |x_0-z|$, so that
$$
\|I^{\a(\cdot)} f_2\|_{L^{q(\cdot)}(\Om\backslash \widetilde{B}(x_0,2t))} \le C \int\limits_{\widetilde{B}(x_0,t)}|f(y)| dy \, \left\||x_0-z|^{\a(z)-n} \right\|_{L^{q(\cdot)}(\Om\backslash \widetilde{B}(x_0,2t))}.
$$
From the log-condition for $\al(\cdot)$ it follows that
$$c_1|x_0-z|^{\a(x_0)-n}\le |x_0-z|^{\a(z)-n} \le c_2 |x_0-z|^{\a(x_0)-n}.$$
Therefore,
\begin{equation}\label{problem}
\|I^{\a(\cdot)} f_2\|_{L^{q(\cdot)}(\Om\backslash \widetilde{B}(x_0,2t))}
\le C \int\limits_{\widetilde{B}(x_0,t)}|f(y)| dy \, \left\||x_0-z|^{\a(x_0)-n} \right\|_{L^{q(\cdot)}(\Om\backslash \widetilde{B}(x_0,2t))}.
\end{equation}
The norm in the integral on the right-hand side is estimated by means of Lemma \ref{lemma},
which yields

$$\|I^{\a(\cdot)} f_2\|_{L^{q(\cdot)}(\Om\backslash \widetilde{B}(x_0,2t))}\le C t^{-\frac{n}{p^\prime(x_0)}} \int\limits_{\widetilde{B}(x_0,t)}|f(y)| dy.
$$
It remains to make use  of \eqref{Ga14} and obtain
\begin{equation}\label{Ha2}
\|I^{\a(\cdot)} f_2\|_{L^{q(\cdot)}(\Om\backslash \widetilde{B}(x_0,2t))}\le C t^{-\frac{n}{p^\prime(x_0)}}\int_0^{t} s^{\frac{n}{p^\prime(x_0)}-1}\|f\|_{L^{p(\cdot)}(\Om\backslash \widetilde{B}(x_0,s))}ds.
\end{equation}

From \eqref{Ha1} and \eqref{Ha2} we arrive at \eqref{GOPRiesz}.
\end{proof}

\begin{theorem}\label{M1X}
Let $\Om \subset \Rn$ be an open bounded  set, $x_0 \in \Om$ and  $p(\cdot), \al(\cdot) \in WL(\Om)$ satisfy assumption \eqref{h0} and   \eqref{1}, $q(x)$ given by \eqref{1cxz}  and the functions $\om_1(r)$ and $\om_2(r)$ fulfill  the condition
\begin{equation}\label{eq3.6.VZX}
t^{\a(x_0)} \int_0^{t}  \om_1(r)\frac{dr}{r} \le C \,\om_2(t),
\end{equation}
where $C$ does not depend on $t$. Then the operators $M^{\al(\cdot)}$
and $I^{\al(\cdot)}$ are bounded from ${\dual \cal M}_{\{x_0\}}^{p(\cdot),\om_1}(\Om)$ to
${\dual \cal M}_{\{x_0\}}^{q(\cdot),\om_2}(\Om)$.
\end{theorem}
\begin{proof}
It suffices to prove the boundedness of the operator  $I^{\al(\cdot)}$, since $M^{\al(\cdot)} f(x)\le C I^{\al(\cdot)}|f|(x)$.

Let $f\in {\dual \cal M}_{\{x_0\}}^{p(\cdot),\om}(\Om)$.
 We have
 \begin{equation}\label{ggyxx}
\|I^{\al(\cdot)} f\|_{{\dual \cal M}^{\{x_0\}}_{q(\cdot),\om_2}(\Om)}=\sup_{t>0}
\frac{t^{\frac{n}{q^\prime(x_0)}}}{\om_2(t)} \|\chi_{\Om\backslash \widetilde{B}(x_0,t)} I^{\al(\cdot)}
f\|_{L^{q(\cdot)}(\Om)}.
\end{equation}
We estimate $\|\chi_{\Om\backslash \widetilde{B}(x_0,t)}I^{\al(\cdot)}
f\|_{L^{q(\cdot)}(\Om)}$
in \eqref{ggyxx} by means of Theorem \ref{G4}. This theorem is applicable since the integral \eqref{converg} with $\om=\om_1$ converges by
 \eqref{qwet}. We obtain
\begin{align*}
\|I^{\al(\cdot)} f\|_{{\dual \cal M}^{\{x_0\}}_{q(\cdot),\om_2}(\Om)} & \le
C \sup_{t>0} \frac{t^{-\frac{n}{p^\prime(x_0)}+\frac{n}{q^\prime(x_0)}}}{\om_2(t)} \int_0^{t}
r^{\frac{n}{p^\prime(x_0)}-1}\|f\|_{L^{p(\cdot)}(\Om\backslash \widetilde{B}(x_0,r))} dr
\\
&\le C\|f\|_{{\dual \cal M}^{\{x_0\}}_{p(\cdot),\om_1}(\Om)} \sup_{t>0}
\frac{t^{\a(x_0)}}{\om_2(t)} \int_0^{t}\frac{\om_1(r)}{r}dr.
\end{align*}
It remains to make use of the condition \eqref{eq3.6.VZX}.
\end{proof}
\begin{corollary} \label{VHS2}
 Let $\Om \subset \Rn$ be an open bounded  set and  $p(\cdot), \al(\cdot) \in WL(\Om)$ satisfy
assumption \eqref{h0} and   \eqref{1}, $q(x)$ given by \eqref{1cxz}, $x_0 \in \Om$ and  $\frac{\lambda}{p^\prime(x_0)} \le \frac{\mu}{q^\prime(x_0)}$.
Then the operators $M^{\al(\cdot)}$ and $I^{\al(\cdot)}$ are bounded from
${\dual \cal L}_{\{x_0\}}^{p(\cdot),\lambda}(\Om)$
to ${\dual \cal L}_{\{x_0\}}^{q(\cdot),\mu}(\Om)$.
\end{corollary}

\

\section{Singular integral operators in the spaces
 ${\dual \cal M}_{\{x_0\}}^{p(\cdot),\om}(\Om)$}\label{singular}
\setcounter{theorem}{0} \setcounter{equation}{0}

\

Theorems \ref{HG11} and \ref{SIO1} proved below, in the case of  the constant
exponent $p$  were proved in \cite{Gul},  Theorems 1.4.2 and 1.4.6 (see also \cite{GulBook}, p. 132, 135).

\begin{theorem} \label{HG11}
Let $\Om$ be an open bounded set,  $p\in WL(\Om)$ satisfy
condition \eqref{h0} and $f\in L^{p(\cdot)}(\Om\backslash \widetilde{B}(x_0,t))$ for every $t\in (0,\ell)$.
If the integral
$$
\int_0^\ell
r^{\frac{n}{p^\prime(x_0)}-1}\|f\|_{L^{p(\cdot)}(\Om\backslash \widetilde{B}(x_0,r))} dr
$$
converges, then
\begin{equation*}\label{GOPSI}
\|T f\|_{L^{p(\cdot)}(\Om\backslash \widetilde{B}(x_0,t))} \le
Ct ^{-\frac{n}{p^\prime(x_0)}}\int_{0}^{2t}
r^{\frac{n}{p^\prime(x_0)}-1}\|f\|_{L^{p(\cdot)}(\Om\backslash \widetilde{B}(x_0,r))} dr,
\end{equation*}
where $C$ does not depend on
 $f$, $x_0$ and $t\in (0,\ell)$.
\end{theorem}

\begin{proof} We split the function $f$ in the form  $f_1+f_2$ as as in \eqref{repr} and have
$$
\|Tf\|_{L^{p(\cdot)}(\Om\backslash \widetilde{B}(x_0,2t))}
\le \|Tf_1\|_{L^{p(\cdot)}(\Om\backslash \widetilde{B}(x_0,2t))}
+\|Tf_2\|_{L^{p(\cdot)}(\Om\backslash \widetilde{B}(x_0,2t))}.
$$

Taking into account that $f_1\in L^{p(\cdot)}(\Om)$, by  Theorem \ref{SIO} we have
\begin{equation*}
\|Tf_1\|_{L^{p(\cdot)}(\Om\backslash \widetilde{B}(x_0,2t))}\le
\|Tf_1\|_{L^{p(\cdot)}(\Om)}
 \le C \|f_1\|_{L^{p(\cdot)}(\Om)}=C \|f\|_{L^{p(\cdot)}(\Om\backslash \widetilde{B}(x_0,t))}.
\end{equation*}
Then in view of \eqref{ADI}
\begin{equation}\label{Ga10'}
\|Tf_1\|_{L^{p(\cdot)}(\Om\backslash \widetilde{B}(x_0,t))}\le Ct^{-\frac{n}{p^\prime(x_0)}}
\int_0^{t} r^{\frac{n}{p^\prime(x_0)}-1}\|f\|_{L^{p(\cdot)}(\Om\backslash \widetilde{B}(x_0,r))} dr.
\end{equation}

To estimate $\|Tf_2\|_{L^{p(\cdot)}(\Om\backslash \widetilde{B}(x_0,2t))},$ note  that  $\frac{1}{2} |x_0-z| \le |z-y|\le\frac{3}{2} |x_0-z|$ for $z\in \Om\backslash \widetilde{B}(x_0,2t)$ and $y\in \widetilde{B}(x_0,t)$, so that
\begin{align*}
\|Tf_2\|_{L^{p(\cdot)}(\Om\backslash \widetilde{B}(x_0,2t))} & \le C
\left\|\int_{\widetilde{B}(x_0,t)}|z-y|^{-n}f(y) dy\right\|_{L^{p(\cdot)}(\Om\backslash \widetilde{B}(x_0,2t))}
\\
&\le C \int_{\widetilde{B}(x_0,t)}|f(y)| dy\||x_0-z|^{-n}\|_{L^{p(\cdot)}(\Om\backslash \widetilde{B}(x_0,2t))}.
\end{align*}

Therefore, with the aid of the estimate \eqref{estikmate} and inequality  \eqref{Ga14}, we get
\begin{equation*}\label{Ga11'}
\|Tf_2\|_{L^{p(\cdot)}(\Om\backslash \widetilde{B}(x_0,2t))}\le Ct^{-\frac{n}{p^\prime(x_0)}}\int_0^{t} s^{\frac{n}{p^\prime(x_0)}-1}\|f\|_{L^{p(\cdot)}(\Om\backslash\widetilde{B}(x_0,s))}ds,
\end{equation*}
which together with \eqref{Ga10'} yields \eqref{GOPSI}.
\end{proof}

\begin{theorem}\label{SIO1}
Let $\Om \subset \Rn$ be an open bounded  set, $x_0 \in \Om$,  $p\in WL(\Om)$ satisfy
condition \eqref{h0} and $\om_1(t)$ and $\om_2(t)$ fulfill condition
\eqref{eq3.6.VZ}. Then the singular integral operator $T$ is bounded from the
space ${\dual \cal M}_{\{x_0\}}^{p(\cdot),\om_1} (\Om)$ to the space ${\dual \cal M} ^{\{x_0\}}_{p(\cdot),\om_2}
(\Om)$.
\end{theorem}

\begin{proof} Let $f\in {\dual \cal M}^{\{x_0\}}_{p(\cdot),\om_1}(\Om)$. We follow the procedure already used in the proof of Theorems \ref{M1} and \ref{M1X}: in the norm
\begin{equation}\label{ggy}
\|Tf\|_{{\dual \cal M}_{\{x_0\}}^{p(\cdot),\om_2}(\Om)} = \sup_{t>0} \frac{t^{\frac{n}{p^\prime (x_0)}}}{\om_2(t)} \|Tf\chi_{\Om\backslash \widetilde{B}(x_0,t)}\|_{L^{p(\cdot)}(\Om)},
\end{equation}
we estimate $\|Tf\chi_{\Om\backslash \widetilde{B}(x_0,t)}\|_{L^{p(\cdot)}(\Om)}$  by means of Theorem \ref{HG11}  and obtain
\begin{align*}
\|Tf\|_{{\dual \cal M}_{\{x_0\}}^{p(\cdot),\om_2}(\Om)} &
\le C\sup_{t>0} \frac1{\om_2(t)} \int_0^{t}
r^{\frac{n}{p(x_0)}-1}\|f\|_{L^{p(\cdot)}(\Om\backslash \widetilde{B}(x_0,r))} dr
\\
&\le C\|f\|_{{\dual \cal M}^{\{x_0\}}_{p(\cdot),\om_1}(\Om)} \sup_{t>0}
\frac1{\om_2(t)} \int_0^{t}\om_1(r)\frac{dr}{r} \le C\|f\|_{{\dual \cal M}^{\{x_0\}}_{p(\cdot),\om_1}(\Om)}.
\end{align*}
\end{proof}

\begin{corollary} \label{VHS3}
 Let $\Om \subset \Rn$ be an open bounded  set,  $p\in WL(\Om)$ satisfy the
assumption \eqref{h0}, $x_0 \in \Om$ and  $0 \le \lambda < n$, $\lambda \le \mu \le n$.
Then the singular integral operator $T$ is bounded from
${\dual \cal L}_{\{x_0\}}^{p(\cdot),\lambda}(\Om)$
to ${\dual \cal L}_{\{x_0\}}^{p(\cdot),\mu}(\Om)$.
\end{corollary}

\vs{5mm}
{\bf Acknowledgements.}  The research of V. Guliyev and J. Hasanov was
partially supported by
the grant of Science Development Foundation under the President of the
Republic of Azerbaijan project EIF-2010-1(1)-40/06-1. The research of
V. Guliyev and S. Samko was partially supported by the Scientific and
Technological Research Council of Turkey (TUBITAK Project No:
110T695).

\end{document}